\documentclass[12pt]{amsart}
\newtheorem{defn}{Definition}
\newtheorem{lem}{Lemma}

\newtheorem{thm}{Theorem}

\newcommand{\defi}{\stackrel {\Delta}{=}}
\newcommand{\bk}{b(k,{\bf p})}
\newcommand{\bkl}{\triangle_k b(k,{\bf p}_{\hat l})}
\newcommand{\bkm}{\triangle_k b(k,{\bf p}_{\hat m})}
\newcommand{\bklm}{\triangle^2_k b(k,{\bf p}_{\hat l,\hat m})}
\newcommand{\p}{{\bf p}}
\newcommand{\di}{\vec{\bf u}}
\begin{document}

\title{On the Capacity Region of Multiple Access Adder Channel}
\author{Shaohui Zhai}
\address{{\it Home Address:} 806 Cherry Lane, Apt. 208, East Lansing, MI 48823, USA}
\email{szhai@lynx.neu.edu}
\begin{abstract}
We compute the capacity region of the $s$-user Multiple
Access Adder Channel and obtain an explicit description of this region as
only one polyhedron, instead of a convex hull of infinitely many
polyhedrons. We obtain this by proving a conjecture raised by Lindstr\"{o}m
in 1975, which in turn follows from certain convexity results, both continuous
and discrete.  
\end{abstract}
\maketitle
\section {Preliminaries}
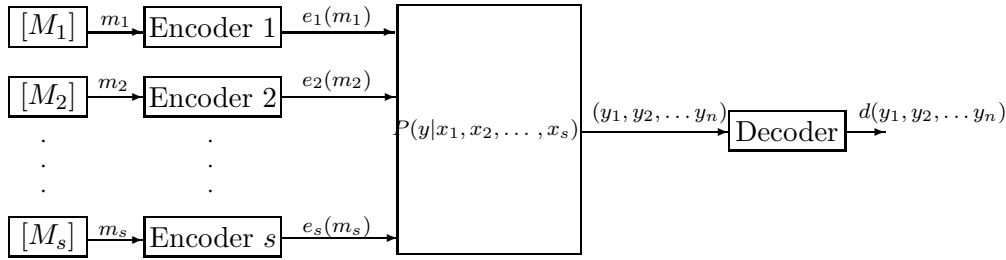
\subsection{Multiple Access Channel}
\begin{defn}
An $s$-user Multiple Access Channel (MAC) consists of $s$ input alphabets, 
${\mathcal X}_1, {\mathcal X}_2, \ldots,{\mathcal X}_s,\;$ an output alphabet 
${\mathcal Y}\;$, and a 
probability transition matrix $P(y\mid x_1, x_2, \ldots, x_s)$, which defines the
probability of receiving the output symbol $y$ given input symbols 
$x_1, x_2, \ldots, x_s$.
\end{defn}
\begin{defn}
An $(M_1,M_2,\ldots,M_s;n)$ code for the $s$-user MAC
consists of $s$ encoding functions $e_1,\ldots,e_s$ which encode each message
$m_i\in [M_i]=\{1,2,\ldots,M_i\}$ into a string of length $n$ in the input 
alphabet $\mathcal{X}_i$:
$$\begin{array}{cccl}
e_i: & [M_i] & \rightarrow & {\mathcal X}^n_i\\
     &  m_i  & \mapsto     & e_i(m_i)=(e^{(1)}_i(m_i), e^{(2)}_i(m_i),
\ldots, e^{(n)}_i(m_i));
\end{array}$$
and a decoding function: 
$$\begin{array}{cccl}
d:\;\;& {\mathcal Y}^n & \rightarrow &[M_1]\times [M_2]\times \cdots \times [M_s]\\
& (y_1,\ldots,y_n) &\mapsto& (\hat{m}_1,\ldots,\hat{m}_s).
\end{array}$$
\end{defn}

The model of communication is as follows (see Fig. 1).  Each user $i$ chooses 
a label $m_i$ 
uniformly from its message set $[M_i]$ and sends the codeword 
$e_i(m_i)$.  All the users send their codewords simultaneously, symbol by 
symbol.  After $n$ transmissions, a sequence of $n$ output symbols $(y_1,y_2,
\ldots,y_n)$ is received according to the transition probability matrix. 
This is then decoded by the decoding function $d$ into a set of messages 
$(\hat{m}_1,\ldots,\hat{m}_s)$.  If the decoded messages are identical 
to the messages sent, the communication has succeeded.  Otherwise,
we say an {\it error} has occurred.\\ 
\vspace{20pt}
\unitlength 0.90mm
\linethickness{0.4pt}
\begin{picture}(153.67,56.99)
\put(28.66,19.99){\framebox(20,5.67)[cc]{\small Encoder $s$}}
\put(28.66,40.66){\framebox(20,5.67)[cc]{\small Encoder $2$}}
\put(28.66,50.99){\framebox(20,5.67)[cc]{\small Encoder $1$}}
\put(66,20.32){\framebox(27,36.67)[cc]
{\tiny $P(y|x_1,x_2,\ldots,x_s)$  }}
\put(115,35.66){\framebox(17.00,5.67)[cc]{\small Decoder}}
\put(8.66,50.99){\framebox(11.33,5.67)[cc]{\small$[M_1]$}}
\put(8.66,40.66){\framebox(11.33,5.67)[cc]{\small$[M_2]$}}
\put(8.66,19.99){\framebox(11.33,5.67)[cc]{\small$[M_s]$}}
\put(13.66,36.99){\circle*{0.00}}
\put(13.66,33.33){\circle*{0.00}}
\put(13.66,29.66){\circle*{0.00}}
\put(38.33,29.66){\circle*{0.00}}
\put(38.33,33.33){\circle*{0.00}}
\put(38.33,36.99){\circle*{0.00}}
\put(24.33,54.66){\makebox(0,0)[cc]{\tiny $ m_1$}}
\put(24.00,45.33){\makebox(0,0)[cc]{\tiny$ m_2$}}
\put(24.00,23.66){\makebox(0,0)[cc]{\tiny$m_s$}}
\put(57.00,55.33){\makebox(0,0)[cc]{\tiny${e_1(m_1)}$}}
\put(57.00,45.99){\makebox(0,0)[cc]{\tiny$e_2(m_2)$}}
\put(57.00,24.33){\makebox(0,0)[cc]{\tiny$e_s(m_s)$}}
\put(104.67,41.16){\makebox(0,0)[cc]
{\tiny$(y_1,y_2,\ldots y_n)$}}
\put(145,41.16){\makebox(0,0)[cc]
{\tiny$d(y_1,y_2,\ldots y_n)$}}
\put(20.00,53.03){\vector(1,0){8.33}}
\put(20.00,43.53){\vector(1,0){8.33}}
\put(20.00,22.33){\vector(1,0){8.33}}
\put(49.00,53.03){\vector(1,0){16.5}}
\put(48.66,43.53){\vector(1,0){17}}
\put(49.00,22.66){\vector(1,0){17}}
\put(132.33,38.33){\vector(1,0){5.67}}
\put(79.00,6.00){\makebox(0,0)[cc]
{{\small Fig. 1. \,Multiple Access Channel}}}
\put(93.33,38.33){\vector(1,0){21.5}}
\end{picture}
\\
\begin{defn}
The average probability of error for an $(M_1,M_2,\ldots,M_s; n)$ code is
$$
P^{(n)}_e=\frac{1}{M_1M_2\cdots M_s}\!\!\!\!\sum_{\tiny \begin{array}{c}
(m_1,\ldots,m_s)\\ \in [M_1]\times \cdots\times[M_s]\end{array}}\!\!\!\!\!
{P_r(\,d(Y^n)\neq (m_1,\ldots,m_s)\,|\,(m_1,\ldots,m_s)\;\mbox {sent})}.
$$\end{defn}
We use {\em rate} to measure the efficiency of an $(M_1,\ldots,M_s;n)$ code.
\begin{defn}
The {\it rate} tuple $(R_1,\ldots,R_s)$ of an 
$(M_1,\ldots,M_s;n)$ code is defined as $$
R_i=\frac{\log_2 M_i}{n}\;\;\;\;\;\mbox {bits per transmission in the i-th 
sender}.$$ 
\end{defn}
\begin{defn} A rate tuple $(R_1,R_2,\ldots,R_s)$ is achievable if there exists 
a sequence of $(2^{nR_1},2^{nR_2},\ldots,2^{nR_s};n)$ codes with 
$P^{(n)}_e\rightarrow 0$ as $n\rightarrow \infty.$ \end{defn}
\begin{defn} The capacity region of the MAC is the closure of the set of all the 
achievable rate tuples.\end{defn}

\begin{defn}
A Multiple Access Channel is called Multiple Access Adder Channel (MAAC) when the 
input alphabets are ${\mathcal X}_1={\mathcal X}_2= \ldots={\mathcal X}_s=\{0,1\}$ 
and the channel output is the sum of the $s$ input symbols:
$Y=X_1+X_2+\cdots +X_s$, so that the output alphabet is
${\mathcal Y}=\{0,1,2,\ldots,s\}$.\end{defn}  
This is a deterministic channel in that a given
s-tuple of inputs always yields the same output. \medskip

\subsection{Capacity Region}
A simple characterization of the capacity region for the two-user MAC 
was first presented by Ahlswede \cite{Ahls1} in 1973.  
Ulrey \cite{Ulrey} generalized it to the MAC with more than two users. 

Let $A\subseteq [s]=\{1,2,\ldots,s\}$, $A^c$ denote the complement of $A$ 
in $[s]$, $R(A)\defi \sum_{i\in A}R_i, X(A)\defi \{X_i: i\in A\},$ and 
$I(X;Y|Z)$ the conditional mutual information of random variable $X$ and $Y$
given random variable $Z$.
\begin{thm}[Ulrey] \label{capacity} The capacity region of the $s$-user
Multiple Access Channel is the closure of the convex hull of the rate tuples 
$(R_1,R_2,\ldots,R_s)$ satisfying
$$R(A)\leq I(X(A);Y|X(A^c))\;\;\;\;\mbox {for all the subsets}\;\;\;
A\subseteq [s]$$
for some product distribution $p_1(x_1)p_2(x_2)\ldots p_s(x_s)$.\end{thm}\medskip

\section{Results}
For the $s$-user MAAC, we now present a finite description of its capacity region.

Let ${\mathcal B}(n, p)$ be the binomial distribution with $n$ independent trials
and probability $p$ of a success in each trial, i.e., if $X\sim {\mathcal B}(n, p)$,
then $Pr(X=k)=\tbinom {n}{k}p^k(1-p)^{n-k}.$ Let $H_n(p)$ denote the Shannon 
entropy of ${\mathcal B}(n, p)$, i.e.,
$$ H_n(p)\defi -\sum^{n}_{k=0}\tbinom {n}{k}p^k(1-p)^{n-k}\cdot 
\log \left(\tbinom {n}{k}p^k(1\!-\!p)^{n-k}\right).$$ 

Now we state our main theorem.
\begin{thm}\label{myconj} The capacity region for the $s$-user Multiple 
Access Adder Channel is the set
$$\{(R_1,R_2,\ldots,R_s):\; R(A)\leq H_{|A|}(\frac 12)\;\;\;\;\mbox{for all}\;\;\;\;A
\subseteq [s]\}.$$\end{thm}
This region is only one 
polyhedron. For example, the capacity region for the 2-user adder 
channel is $$\{(R_1,R_2):\;0\leq R_1\leq 1; 0\leq R_2\leq 1; 
0\leq R_1+R_2\leq 1/2\};$$
the capacity region for the 3-user adder channel is
$$\begin{array}{cl}
\{(R_1,R_2,R_3):& \;R_1\leq 1; R_2\leq 1; R_3\leq 1; R_1+R_2\leq 1/2; 
R_1+R_3\leq 1/2;\\
                & R_2+R_3\leq 1/2; R_1+R_2+R_3\leq (3-\frac 34\log3)\}.
\end{array}$$

A key part of our proof of Theorem~\ref{myconj} is the following result,
conjectured by Lindstr\"{o}m \cite{Lind}:
\begin{thm}\label{lconj} Let $X_1,X_2,\ldots,X_n$ be independent random variables 
$=0\,$ or $1\,$ with $Pr(X_i=1)=p_i\,$ and $\,Pr(X_i=0)=1-p_i,\,$ where 
$\,0\leq p_i\leq 1\,$ for $\,i=1,2,\ldots,n.\,$ The Shannon entropy of the 
random variable 
$X_1+X_2+\cdots+X_n\,$ is then maximal when $\,p_i=1/2\,$ for all 
$\,i=1,2,\ldots,n,\;\;$ i.e., 
$$\max_{0\leq p_1,p_2,\ldots,p_n\leq 1}H(\sum_{i=1}^{n}X_i)=H_n(1/2)$$
\end{thm}

{\bf Remark.} {\it The result in Theorem~\ref{myconj} has apparently obtained by 
Liao \cite{Liao} in his Ph.D. dissertation in 1972, but his proof is never 
published. Our proof will include proving Theorem~\ref{lconj}, conjectured 
in 1976.}\medskip

\section{Proof of Results}
\subsection{Theorem~\ref{lconj} implies Theorem~\ref{myconj}}
To compute the capacity region for the $s$-user adder channel, by 
Theorem~\ref{capacity}, we need to maximize$$I(X(A);Y|X(A^c))$$
for any $A\subseteq [s].\,$

Since $Y=X_1+X_2+\cdots +X_s$ is determined by 
$X(A)$ and $X(A^c)$ for all $A\subseteq [s]\,$, the conditional Shannon entropy
of $Y$ given $X(A)$ and $X(A^c)$ is zero, i.e.,  
$H\left(Y|X(A),X(A^c)\right)=0$, so
$$\begin{array}{ccl}
\displaystyle I(X(A);Y|X(A^c))&=&
\displaystyle H\left(Y|X(A^c)\right)-H\left(Y|X(A),X(A^c)\right)\\
&=&\displaystyle H\left(Y|X(A^c)\right)\\
&=&\displaystyle \sum_{x}{H\left(Y|X(A^c)\!=\!x\right)
                   \cdot P_r\left(X(A^c)\!=\!x\right)}\\
&=&\displaystyle \sum_{x}{H(\sum_{i\in A} X_i\!+\!\sum_{i \in A^c} 
                X_i\;|\;X(A^c)\!=\!x)\cdot P_r\left(X(A^c)=x\right)}\\
&=&\displaystyle \sum_{x} {H(\sum_{i\in A}X_i)\cdot 
                      P_r\left(X(A^c)\!=\!x\right)}\\
&=&\displaystyle H(\sum_{i\in A}X_i)\cdot \sum_{x}
                      {P_r\left(X(A^c)\!=\!x\right)}\\
&=&\displaystyle H(\sum_{i\in A}X_i).
\end{array}$$
So, $$\max_{0\leq p_1,p_2,\ldots,p_s\leq 1}I(X(A);Y|X(A^c))=
\max_{\{i\in A:\;0\leq p_i\leq 1\}}H(\sum_{i\in A}X_i).$$
Thus, if Theorem~\ref{lconj} is true, then $H\left(\sum_{i\in A}X_i \right)$ is 
maximized when $p_i=1/2,\; i\in A$ for all $A\subseteq [s].\,$ 
But when $p_i=1/2,\; i\in A,\;$ the random variable $\sum_{i\in A}X_i$ 
has the binomial distribution with parameter $(|A|,1/2)$, i.e., 
$$\max_{0\leq p_1,p_2,\ldots,p_s\leq 1}I(X(A);Y|X(A^c))=
\max_{\{i\in A:\;0\leq p_i\leq 1\}}H(\sum_{i\in A}X_i)=H_{|A|}(1/2)$$ 
for all $A\subseteq [s].\,$ This, together with Theorem~\ref{capacity},
implies Theorem~\ref{myconj}.

\subsection{Proof of Theorem~\ref{lconj} (Lindstr\"{o}m's conjecture)}
Let $p_i=P_r(X_i=1), 0\leq p_i\leq 1,$ and $q_i=1-p_i,$ for $i=1,2,\ldots,n$.
Let ${\bf p}=(p_1, p_2,\ldots,p_n)$, ${\bf p}_{\hat i}=(p_1,\ldots,p_{i-1},
p_{i+1},\ldots,p_n),\;$ and more generally, for $0\leq r\leq n,\;$ 
${\bf p}_{\hat i_1,\hat i_2,\ldots,\hat i_r}=
(p_{j_1}, p_{j_2},\ldots,p_{j_{n-r}}),$ where $1\leq j_1<j_2<\cdots<j_{n-r}\leq n$ 
and $\{j_1, j_2,\ldots, j_{n-r}\}=[n]\setminus \{i_1, i_2,\ldots, i_r\}.$
 
We need the following result of Mateev \cite{Mateev}. 

\begin{lem}[Mateev]\label{mateev} The Shannon entropy of the binomial 
distribution ${\mathcal B}(n, p)$ reaches its maximum when $p=1/2$, that is,
$$\max_{0\leq p\leq 1} H_n(p)=H_n(1/2).$$
\end{lem}

Let $Z_{\p}=\displaystyle\sum^{n}_{i=1}X_i$. 
We want to show that the Shannon entropy $H(Z_{\p})$ of $Z_{\p}$ is maximized 
when $p_i=1/2$ for all $i=1,2,\ldots,n$, i.e.,
$$\max_{0\leq p_1,p_2,\ldots,p_n\leq 1}H(Z_{\p})=H_n(1/2).$$
Notice that, when $p_1=p_2=\cdots=p_n=p,\;$ $Z_{\p}$ has the binomial distribution 
with parameter $(n,p)$, and $H(Z_{\p})=H_n(p).$ By Lemma~\ref{mateev},  
$\max_{0\leq p\leq 1} H_n(p)=H_n(1/2).$ Therefore, to prove Theorem~\ref
{lconj}, we only need to show that
\begin{equation}\label{goal}\max_{0\leq p_1,p_2,\ldots,p_n\leq 1}H(Z_{\p})=
\max_{0\leq p_1=p_2=\cdots=p_n\leq 1}H(Z_{\p})\end{equation}
We introduce the following lemma,

\begin{lem}\label{concave} Let ${\bf p'}=(p'_1,\ldots,p'_n)$ be any point in
${\bf R^n}$ with $0\leq p'_1,\ldots,p'_n \leq 1$. Let $1\leq l< m\leq n$, 
and $\di=(0,\ldots,\stackrel{l}{1},0,\ldots,\stackrel{m}
{-1},0,\ldots,0)$. Then the function $H(Z_{\p})$ restricted to the line 
$\p={\bf p'}+ t\cdot \di,\;$ where $t$ is any real number such that 
${\bf 0} \leq {\bf p} \leq {\bf 1},\,$ has a unique maximum at $p_l=p_m$.
\end{lem}
Theorem~\ref{lconj} follows immediately from this lemma (which we will prove in the 
next paragraph). Indeed, let ${\bf p^*}=(p_1^*,p_2^*,\ldots,p_n^*)$ be any maximal 
point of $H(Z_{\p})$. Then for any $1\leq l < m\leq n$, $\;{\bf p^*}$ is also the
maximal point of the function $H(Z_{\p})$ restricted to the line 
$\p={\bf p^*}+ t\cdot \di\;$. 
By Lemma~\ref{concave}, $\;p_l^*=p_m^*\;$, for any $l, m$, which proves
(~\ref{goal}).

{\bf Proof of Lemma~\ref{concave}.} We first notice that $H(Z_{\p})$ is symmetric 
in $p_1,p_2,\\ \ldots,p_n$, so the function $H(Z_{\p})$ restricted to the line 
$\p={\bf p'}+ t\cdot \di$ is symmetric in $p_l$ and $p_m$. We will prove that 
the function $H(Z_{\p})$ restricted to the line $\p={\bf p'}+ t\cdot \di$ is 
concave by showing that the second directional derivative of $H(Z_{\p})$ in the 
direction $\di$ is
strictly negative. The concavity guarantees a unique maximum; 
by the symmetry, this maximum occurs at $p_l=p_m$.

Thus, to prove Lemma~\ref{concave}, and hence Theorem~\ref{lconj} and Theorem~\ref
{myconj}, we must show: $$\begin{array}{c}
D^2_{\di}H(Z_{\p})<0  \\
\forall \;\p=(p_1,p_2,\ldots,p_n)\;\;
                      \mbox {with}\;\;0\leq p_1,p_2,\ldots,p_n\leq 1, 
\end{array}$$
where $D_{\di}=\frac{\partial}{\partial p_l}-\frac{\partial}{\partial p_m}.$

The probability distribution of $Z_{\p}$ is described by \\ $$
\bk\defi Pr(Z_{\p}=k)=\left\{\begin{array}{l}
    \displaystyle \sum_{A\subseteq [n]\atop |A|=k}{\prod_{i\in A}p_i
\prod_{j\in A^c}q_j}\;\;\;\;\;\;\;\mbox {for}\;\;0\leq k\leq n \\
\;\;\;\;0 \;\;\;\;\;\;\;\;\;\;\;\;\;\;\;\;\;\;\;\;\;\;\;\;\;\;\;
\mbox {otherwise},\end{array} \right. $$\\
where $A^c$ is the complement of $A$ in $[n].$ Particularly, 
$$b(0, \p)=\prod_{j=1}^{n}q_j,\;\;\;\;\;\;\;\;
b(n,\p)=\prod_{i=1}^{n}p_i.$$
In general, the probability function 
$Pr(Z_{{\bf p}_{\hat i_1,\hat i_2,\ldots,\hat i_r}}=k)\;$ of 
$\;Z_{{\bf p}_{\hat i_1,\hat i_2,\ldots,\hat i_r}}\;$ ($\defi\sum_{1\leq i\leq n\atop
i\neq i_1,\ldots,i_r}X_i$) is 
$$b(k, {\bf p}_{\hat i_1,\hat i_2,\ldots,\hat i_r}) 
=\left\{\begin{array}{l}
    \displaystyle \sum_{A\subseteq [n]\setminus\{i_1,\ldots,i_r\}\atop |A|=k}
{\prod_{i\in A}p_i\prod_{j\in A^c}q_j}\;\;\;\;\;\;\;
\mbox {for}\;\;0\leq k\leq (n-r) \\
\;\;\;\;\;\;\;\;0 \;\;\;\;\;\;\;\;\;\;\;\;\;\;\;\;\;\;\;\;\;\;\;\;\;\;\;\;\;\;\;\;
\;\;\mbox {otherwise},\end{array} \right. $$
where $A^c$ is the complement of $A$ in $[n]\setminus\{i_1,i_2,\ldots,i_r\}.$

The Shannon entropy of $Z_{\p}$ is
$$ H(Z_{\p})= -\sum_{k=0}^n \bk\log \bk.$$
To compute the derivatives of $H(Z_{\p})$, we will repeatedly use the following 
formulas:\\[-1.5em]
$$\begin{array}{ccl}
\displaystyle \frac{\partial \bk}{\partial p_l} &=& \displaystyle \sum_
{l\in A\subseteq [n]\atop |A|=k}\prod_{i\neq l\atop i\in A}p_i
\prod_{j\in A^c}q_j-\sum_{l\notin B\subseteq [n]\atop |B|=k}
\prod_{i\in B}p_i\prod_{j\neq l\atop j\in B^c}q_j\\
&=& \displaystyle \sum_
{A'\subseteq [n]\setminus\{l\}\atop |A'|=k-1}\prod_{i\in A'}p_i
\prod_{j\in (A')^c}q_j-\sum_{B'\subseteq [n]\setminus\{l\}\atop |B'|=k}
\prod_{i\in B'}p_i\prod_{j\in (B')^c}q_j\\
&=& b(k\!-\!1, {\bf p}_{\hat l})- b(k, {\bf p}_{\hat l})\\
&=& -\bkl,
\end{array}$$
where $\triangle_k f(k)=f(k)-f(k\!-\!1),\;$ the difference operator. 
On the boundaries $k=0, k=n$,
the above formula still holds:\\[-1em]
$$\begin{array}{cl}
\displaystyle\frac{\partial b(0,\p)}{\partial p_l}=&\displaystyle
\!\!\!\!-\!\!\prod_{1\leq j\leq n\atop j\neq l}q_j
=-b(0,{\bf p}_{\hat l})=b(-1,{\bf p}_{\hat l})-b(0,{\bf p}_{\hat l})
=-\triangle_k b(0, {\bf p}_{\hat l})\\[1.5em]
\displaystyle\frac{\partial b(n,\p)}{\partial p_l}=& \displaystyle
\!\!\prod_{1\leq i\leq n\atop i\neq l}p_i
=b(n-1,{\bf p}_{\hat l})=b(n-1,{\bf p}_{\hat l})-b(n,{\bf p}_{\hat l})
=-\triangle_k b(n, {\bf p}_{\hat l}).
\end{array}$$
Let us compute the second derivatives of $b(k,\p)$:
$$\begin{array}{rrl}
\displaystyle \frac{\partial^2 \bk}{\partial p_l\partial p_m}&=&
\displaystyle\frac{\partial^2 \bk}{\partial p_m\partial p_l}= 
\frac{\partial}{\partial p_m}(\;b(k\!-\!1,{\bf p}_{\hat l})- 
b(k,{\bf p}_{\hat l})\;)\\
&=&\displaystyle 
   b(k\!-\!2, {\bf p}_{\hat l,\hat m})- b(k\!-\!1, {\bf p}_{\hat l,\hat m})
   -b(k\!-\!1, {\bf p}_{\hat l,\hat m})+ b(k, {\bf p}_{\hat l,\hat m})\\
&=& \bklm,\\[.8em]
\displaystyle \frac{\partial^2 \bk}{\partial p_l\partial p_l}&=& 
\displaystyle\frac{\partial}{\partial p_l}(\;b(k\!-\!1,{\bf p}_{\hat l})- 
b(k,{\bf p}_{\hat l})\;)=0.
\end{array}$$
When $k=0$ and $k=n$, we have $$\begin{array}{ccl}
\displaystyle\frac{\partial^2 b(0,\p)}{\partial p_m\partial p_l}&=&\displaystyle 
\frac{\partial}{\partial p_m}(\;b(-1,{\bf p}_{\hat l})-b(0,{\bf p}_{\hat l})\;)
=\frac{\partial}{\partial p_m}(\;0-\prod_{1\leq j\leq n\atop j\neq l}q_j)\\
&=&\displaystyle\prod_{1\leq j\leq n\atop j\neq l,m}q_j= 
   b(0, {\bf p}_{\hat l,\hat m})\\
&=& b(-2,{\bf p}_{\hat l,\hat m})- 2b(-1, {\bf p}_{\hat l,\hat m})
   + b(0, {\bf p}_{\hat l,\hat m})
   =\triangle_k^2 b(0,{\bf p}_{\hat l,\hat m}),\\[1em]
\displaystyle\frac{\partial^2 b(n,\p)}{\partial p_m\partial p_l}&=&\displaystyle 
\frac{\partial}{\partial p_m}(\;b(n-1,{\bf p}_{\hat l})-b(n,{\bf p}_{\hat l})\;)
=\frac{\partial}{\partial p_m}(\;\prod_{1\leq i\leq n\atop i\neq l}p_i-0\;)\\
&=&\displaystyle\prod_{1\leq i\leq n\atop i\neq l,m}p_i= 
   b(n-2, {\bf p}_{\hat l,\hat m})\\
&=& b(n-2,{\bf p}_{\hat l,\hat m})- 2b(n-1, {\bf p}_{\hat l,\hat m})
   + b(n, {\bf p}_{\hat l,\hat m})
   =\triangle_k^2 b(n,{\bf p}_{\hat l,\hat m}).
\end{array}$$
Then let us compute the derivatives of $H(Z_{\p})$:
$$\begin{array}{lll}
 \displaystyle \frac{\partial H(Z_{\p})}{\partial p_l}&=& \displaystyle -
\sum_{k=0}^n \frac{\partial \bk}{\partial p_l}\cdot \log \bk -\sum_{k=0}^n 
\bk\cdot\frac{1}{\bk}\cdot\frac{\partial \bk}{\partial p_l}\cdot\frac{1}{\ln 2}
\\[2em]
&=&\displaystyle -\sum_{k=0}^n \frac{\partial \bk}{\partial p_l}\cdot \log \bk-
\frac{1}{\ln 2}\frac{\partial}{\partial p_l}\sum_{k=0}^n \bk\\[2em]
&=& \displaystyle -\sum_{k=0}^n \frac{\partial \bk}{\partial p_l}\cdot \log \bk-
\frac{1}{\ln 2}\frac{\partial}{\partial p_l}(1)\\[2em]
&=& \displaystyle -\sum_{k=0}^n \frac{\partial \bk}{\partial p_l}\cdot \log \bk\\[1em]
&=&\displaystyle\sum_{k=0}^n \bkl \cdot \log \bk;
\end{array}$$
and the second derivatives of $H(Z_{\p})$:
$$\begin{array}{lll}
\displaystyle \frac{\partial^2 H(Z_{\p})}{\partial p_l\partial p_m}& = & 
\displaystyle \!\!-\sum_{k=0}^{n}\!\frac{\partial^2 \bk}{\partial p_l\partial p_m}
\!\cdot\!\log \bk-\sum_{k=0}^{n}\!\frac{\partial \bk}{\partial p_l}\!\cdot\!
\frac{1}{\ln 2 \!\cdot\! \bk}\cdot \frac{\partial \bk}{\partial p_m}\\[2em]
&=&\displaystyle \!\!-\sum_{k=0}^{n}\!\bklm\cdot \log \bk-\sum_{k=0}^{n}\!\frac
{\bkl \cdot \bkm }{\ln 2 \cdot \bk},\\[2em]
\displaystyle \frac{\partial^2 H(Z_{\p})}{\partial p_l\partial p_l}& = &
\displaystyle -\sum_{k=0}^{n}\frac{(\bkl)^2}{\ln 2 \cdot \bk}.
\end{array}$$\\[-.8em]
The second directional derivative of $H(Z_{\p})$ in the direction $\di$ is
$$\begin{array}{cl}
&\displaystyle D^2_{\di}H(Z_{\p})
=\displaystyle \frac{\partial^2 H(Z_{\p})}{\partial p_l^2}-2\frac{\partial^2 H(Z_{\p})}
{\partial p_l\partial p_m}+\frac{\partial^2 H(Z_{\p})}{\partial p_m^2}\\[2em]
&=\displaystyle-\sum_{k=0}^{n}\frac{(\bkl)^2}{\ln 2 \cdot \bk}+
  2\sum_{k=0}^{n}\bklm\cdot \log \bk\\[2em]
& \displaystyle \;\;\;\;+2\sum_{k=0}^{n}\frac{\bkl \cdot \bkm }{\ln 2 \cdot \bk}
  -\sum_{k=0}^{n}\frac{(\bkm)^2}{\ln 2 \cdot \bk}\\[2em]
&=\displaystyle \;2\sum_{k=0}^{n}\bklm
  \cdot \log \bk-\sum_{k=0}^{n}\frac{(\bkl-\bkm)^2}{\ln 2\cdot\bk}.
\end{array}$$\\[-.8em]
The second term of $D^2_{\di}$ is clearly negative. We will prove that the first
term is strictly negative, therefore $\;D^2_{\di}H(Z_{\p})<0\;$ will be proved.

To evaluate the first term, we use ``summation by parts":
$$\displaystyle \sum_{k=0}^{n} (\triangle_k f(k))\cdot g(k)=
-\sum_{k=0}^n f(k\!-\!1)\cdot(\triangle_k g(k))\;\;\;\mbox {provided}\;\;f(n)=f(-1)=0.
$$ 
Noticing that $$\begin{array}{c}
\displaystyle b(k,{\bf p}_{\hat l,\hat m})=0 \;\;\;\;\;\mbox{when}\;\; k<0\;\;
\mbox{or}\;\;k>n-2\\
\triangle_k b(k,{\bf p}_{\hat l,\hat m})=0\;\;\;\;\mbox{when}\;\;k<0\;\;\mbox{or}
\;\;k>n-1,
\end{array}$$ 
we can apply summation by parts twice to the first term:
$$\begin{array}{cl}
&\displaystyle\sum_{k=0}^{n}\bklm \cdot \log \bk\\[1.5em]
=&\displaystyle-\displaystyle \sum_{k=0}^n \triangle_k b(k\!-\!1,
{\bf p}_{\hat l,\hat m})\cdot\triangle_k(\,\log \bk \,)\\[1.5em]
=&\displaystyle \sum_{k=0}^n b(k\!-\!2,{\bf p}_{\hat l,\hat m})\cdot 
  \triangle^2_k (\,\log \bk\,)
\end{array}$$\\[-.8em]
Since $b(k,{\bf p}_{\hat l,\hat m})$'s are probabilities, they are 
non-negative, and at least one of them is positive. Finally, we have reduced our
proof to the following lemma:
\begin{lem}\label{logcon}
$\bk$ is log-concave as a function of k, that is, $$
\triangle^2_k (\,\log \bk\,) < 0,\;\;\;\;\mbox{for any}\;\;\; {\bf 0}\leq \p\leq
{\bf 1}.$$
\end{lem}

{\bf Proof of Lemma~\ref{logcon}.} Since
$$\begin{array}{lll}
\displaystyle \triangle^2_k (\,\log \bk\,) &=& \log \bk-2\log b(k\!-\!1,{\bf p})+
                    \log b(k\!-\!2,{\bf p})\\[1em]
&=& \displaystyle \log \frac{b(k\!-\!2, {\bf p})\cdot b(k,{\bf p})}
{(b(k\!-\!1,{\bf p}))^2},
\end{array}$$
we only need to show that 
$$b(k\!-\!2,{\bf p})\cdot b(k,{\bf p}) < (b(k\!-\!1,{\bf p}))^2,$$ 
or equivalently, 
$$\;b(k\!-\!1,{\bf p})\cdot b(k\!+\!1,{\bf p}) < (b(k,{\bf p}))^2.$$ 
Let us expand the left-hand side of the inequality:
$$\begin{array}{cl}
& b(k\!-\!1,{\bf p})\cdot b(k\!+\!1,{\bf p})\\[1em]
=&\!\!\displaystyle 
\sum_{A\subseteq [n]\atop |A|=k\!-\!1}{\prod_{a\in A}p_a\prod_{a'\in A^c}q_{a'}}
\cdot
\sum_{B\subseteq [n]\atop |B|=k\!+\!1}{\prod_{b\in B}p_b\prod_{b'\in B^c}q_{b'}}
\\[1.5em]
=&\!\! \displaystyle\sum_{A,B\subseteq [n]
\atop |A|=k\!-\!1,|B|=k\!+\!1}\prod_{a\in A}p_a\prod_{b\in B}p_b\prod_{a'\in A^c}
q_{a'}\prod_{b'\in B^c}q_{b'}\\[1.5em]
=&\!\!\displaystyle\sum_{A,B\subseteq [n]\atop |A|=k\!-\!1,|B|=k\!+\!1}
\prod_{A\cap B}p^2_i
\prod_{(A\cup B)\setminus (A\cap B)}p_{i'}\prod_{(A\cup B)^c}q^2_j\prod
_{(A\cap B)^c\setminus (A\cup B)^c}q_{j'}\\[2em] 
\stackrel{(1)}{=}&\!\!\displaystyle\sum_{w=0}^{k-1}\sum_{I\subseteq [n]\atop |I|=w}\sum_
{J\subseteq [n]\setminus I\atop |J|=2k\!-\!2w}\!\!\!\tbinom {2k-2w}{k-w-1}
\prod_{I}p^2_i\prod_{J}p_{i'}
\prod_{(I\cup J)^c}q_j^2\prod_{I^c\setminus (I\cup J)^c}q_{j'}.
\end{array}$$\\[-1em]
Equality $(1)$ may be illustrated by Figure 1.\\[2.5em]
\unitlength 1mm
\linethickness{0.4pt}
\begin{picture}(120.00,90.00)
\put(30.00,37.00){\framebox(80.00,55.00)}
\multiput(60.00,80.06)(0.99,-0.10){3}{\line(1,0){0.99}}
\multiput(62.98,79.76)(0.36,-0.11){8}{\line(1,0){0.36}}
\multiput(65.84,78.88)(0.22,-0.12){12}{\line(1,0){0.22}}
\multiput(68.48,77.45)(0.14,-0.11){17}{\line(1,0){0.14}}
\multiput(70.77,75.52)(0.12,-0.15){16}{\line(0,-1){0.15}}
\multiput(72.64,73.18)(0.11,-0.22){12}{\line(0,-1){0.22}}
\multiput(74.01,70.52)(0.12,-0.41){7}{\line(0,-1){0.41}}
\multiput(74.83,67.63)(0.11,-1.49){2}{\line(0,-1){1.49}}
\multiput(75.06,64.65)(-0.09,-0.74){4}{\line(0,-1){0.74}}
\multiput(74.69,61.67)(-0.12,-0.36){8}{\line(0,-1){0.36}}
\multiput(73.74,58.83)(-0.11,-0.20){13}{\line(0,-1){0.20}}
\multiput(72.25,56.23)(-0.12,-0.13){17}{\line(0,-1){0.13}}
\multiput(70.27,53.98)(-0.15,-0.11){16}{\line(-1,0){0.15}}
\multiput(67.88,52.17)(-0.25,-0.12){11}{\line(-1,0){0.25}}
\multiput(65.19,50.86)(-0.41,-0.11){7}{\line(-1,0){0.41}}
\multiput(62.29,50.12)(-1.50,-0.08){2}{\line(-1,0){1.50}}
\multiput(59.29,49.96)(-0.74,0.11){4}{\line(-1,0){0.74}}
\multiput(56.33,50.39)(-0.31,0.11){9}{\line(-1,0){0.31}}
\multiput(53.51,51.41)(-0.20,0.12){13}{\line(-1,0){0.20}}
\multiput(50.95,52.96)(-0.13,0.12){17}{\line(-1,0){0.13}}
\multiput(48.75,54.99)(-0.12,0.16){15}{\line(0,1){0.16}}
\multiput(46.99,57.42)(-0.11,0.25){11}{\line(0,1){0.25}}
\multiput(45.74,60.15)(-0.11,0.49){6}{\line(0,1){0.49}}
\put(45.07,63.06){\line(0,1){2.99}}
\multiput(44.98,66.06)(0.10,0.59){5}{\line(0,1){0.59}}
\multiput(45.48,69.01)(0.11,0.28){10}{\line(0,1){0.28}}
\multiput(46.57,71.80)(0.12,0.18){14}{\line(0,1){0.18}}
\multiput(48.18,74.33)(0.12,0.12){18}{\line(0,1){0.12}}
\multiput(50.26,76.49)(0.16,0.11){15}{\line(1,0){0.16}}
\multiput(52.73,78.19)(0.28,0.12){10}{\line(1,0){0.28}}
\multiput(55.48,79.36)(0.75,0.12){6}{\line(1,0){0.75}}
\multiput(80.00,77.87)(1.13,-0.11){3}{\line(1,0){1.13}}
\multiput(83.39,77.54)(0.36,-0.11){9}{\line(1,0){0.36}}
\multiput(86.65,76.58)(0.22,-0.11){14}{\line(1,0){0.22}}
\multiput(89.68,75.02)(0.15,-0.12){18}{\line(1,0){0.15}}
\multiput(92.35,72.91)(0.12,-0.14){19}{\line(0,-1){0.14}}
\multiput(94.57,70.34)(0.11,-0.20){15}{\line(0,-1){0.20}}
\multiput(96.27,67.38)(0.11,-0.32){10}{\line(0,-1){0.32}}
\multiput(97.37,64.17)(0.12,-0.84){4}{\line(0,-1){0.84}}
\multiput(97.85,60.80)(-0.09,-1.70){2}{\line(0,-1){1.70}}
\multiput(97.68,57.40)(-0.12,-0.47){7}{\line(0,-1){0.47}}
\multiput(96.86,54.09)(-0.12,-0.26){12}{\line(0,-1){0.26}}
\multiput(95.44,51.00)(-0.12,-0.16){17}{\line(0,-1){0.16}}
\multiput(93.45,48.24)(-0.12,-0.12){20}{\line(-1,0){0.12}}
\multiput(90.97,45.90)(-0.18,-0.11){16}{\line(-1,0){0.18}}
\multiput(88.10,44.08)(-0.29,-0.11){11}{\line(-1,0){0.29}}
\multiput(84.93,42.83)(-0.56,-0.10){6}{\line(-1,0){0.56}}
\put(81.59,42.20){\line(-1,0){3.40}}
\multiput(78.19,42.23)(-0.56,0.11){6}{\line(-1,0){0.56}}
\multiput(74.85,42.89)(-0.29,0.12){11}{\line(-1,0){0.29}}
\multiput(71.70,44.18)(-0.18,0.12){16}{\line(-1,0){0.18}}
\multiput(68.85,46.04)(-0.12,0.12){20}{\line(-1,0){0.12}}
\multiput(66.40,48.41)(-0.11,0.16){17}{\line(0,1){0.16}}
\multiput(64.45,51.20)(-0.12,0.26){12}{\line(0,1){0.26}}
\multiput(63.06,54.31)(-0.11,0.47){7}{\line(0,1){0.47}}
\multiput(62.29,57.62)(-0.06,1.70){2}{\line(0,1){1.70}}
\multiput(62.16,61.02)(0.10,0.67){5}{\line(0,1){0.67}}
\multiput(62.68,64.39)(0.11,0.32){10}{\line(0,1){0.32}}
\multiput(63.83,67.59)(0.12,0.20){15}{\line(0,1){0.20}}
\multiput(65.56,70.52)(0.12,0.13){19}{\line(0,1){0.13}}
\multiput(67.82,73.07)(0.15,0.12){18}{\line(1,0){0.15}}
\multiput(70.51,75.14)(0.23,0.12){13}{\line(1,0){0.23}}
\multiput(73.56,76.67)(0.59,0.11){11}{\line(1,0){0.59}}
\put(40,85){\makebox(0,0)[cc]{$[n]$}}
\put(45.67,51.67){\makebox(0,0)[cc]{$A$}}
\put(50,46.00){\makebox(0,0)[cc]{\small$|A|\!=\!k\!-\!1$}}
\put(101.00,52.00){\makebox(0,0)[cc]{$B$}}
\put(102.00,46.33){\makebox(0,0)[cc]{\small$|B|\!=\!k\!+\!1$}}
\put(67.67,64.00){\makebox(0,0)[cc]{$I$}}
\put(68.33,59.00){\makebox(0,0)[cc]{\small$|I|\!=\!w$}}
\put(71.33,85.67){\makebox(0,0)[cc]{$J$}}
\put(86,85.67){\makebox(0,0)[cc]{\small$|J|\!=\!2k\!-\!2w$}}
\put(70, 30){\makebox(0,0)[cc]{Figure 1.}}
\put(53.67,78.33){\line(-3,-5){8.20}}
\put(59.00,79.67){\line(-3,-5){12.40}}
\put(64.00,79.00){\line(-3,-5){14.80}}
\put(68.67,76.67){\line(-3,-5){15.20}}
\put(62.33,55.67){\line(-3,-5){3.40}}
\put(78.00,77.67){\line(-2,-3){3.78}}
\put(83.67,77.33){\line(-3,-5){8.00}}
\put(89.33,75.00){\line(-3,-5){18.20}}
\put(94.67,70.00){\line(-3,-5){16.60}}
\put(97.67,62.67){\line(-3,-5){11.80}}
\put(69.33,53.00){\line(-2,-3){2.67}}
\put(68.33,85.00){\vector(-1,-1){9.33}}
\put(72.67,83.67){\vector(3,-4){8.00}}
\end{picture}

The right-hand side of the inequality can be simplified as:  
$$\begin{array}{cl}
& (\,b(k,{\bf p})\,)^2\\
=&\!\!\displaystyle 
\sum_{C\subseteq [n]\atop |C|=k}{\prod_{c\in C}p_{c}
\prod_{c'\in C^c}q_{c'}}\cdot
\sum_{D\subseteq [n]\atop |D|=k}{\prod_{d\in D}p_{d}
\prod_{d'\in D^c}q_{d'}}\\
=&\!\! \displaystyle\sum_{C,D\subseteq [n]
\atop |C|=|D|=k}\prod_{c\in C}p_c\prod_{d\in D}p_d\prod_{c'\in C^c}q_{c'}
\prod_{d'\in D^c}q_{d'}\\
=& \displaystyle\sum_{C=D\subseteq [n]\atop |C|=k}\prod_{C}p_i^2
\prod_{C^c}q_j^2 \;\;+\sum_{C\neq D\subseteq [n]\atop |C|=|D=k}\prod_{C}p_c
\prod_{D}p_d\prod_{C^c}q_{c'}\prod_{D^c}q_{d'}\\
\geq & \displaystyle\sum_{C\neq D\subseteq [n]\atop |C|=|D=k}\prod_{C}p_c
\prod_{D}p_d\prod_{C^c}q_{c'}\prod_{D^c}q_{d'} \\[1.8em]
\stackrel{(2)}{=}& \displaystyle\sum_{w=0}^{k-1}\sum_{I\subseteq [n]\atop |I|=w}\sum_
{J\subseteq [n]\setminus I\atop |J|=2k\!-\!2w}\!\!\!\tbinom {2k-2w}{k-w}
\prod_{I}p^2_i\prod_{J}p_{i'}
\prod_{(I\cup J)^c}q_j^2\prod_{I^c\setminus (I\cup J)^c}q_{j'}\\[1.8em]
\stackrel{(3)}{>} &\displaystyle\sum_{w=0}^{k-1}\sum_{I\subseteq [n]\atop |I|=w}\sum_
{J\subseteq [n]\setminus I\atop |J|=2k\!-\!2w}\!\!\!\tbinom {2k-2w}{k-w-1}
\prod_{I}p^2_i\prod_{J}p_{i'}
\prod_{(I\cup J)^c}q_j^2\prod_{I^c\setminus (I\cup J)^c}q_{j'}\\
=&\displaystyle b(k\!-\!1,{\bf p})\cdot b(k\!+\!1,{\bf p}). 
\end{array}$$
Equality $(2)$ can be illustrated by Figure 2., and inequality $(3)$ is true 
since $\binom {2n}{n}>\binom {2n}{n-1}.$ Therefore Lemma~\ref{logcon} is proved.\\
[6em]
\unitlength 1mm
\linethickness{0.4pt}
\begin{picture}(100.00,43.50)
\put(30,9){\framebox(80.00,55.00)}
\multiput(60.00,52.06)(0.99,-0.10){3}{\line(1,0){0.99}}
\multiput(62.98,51.76)(0.36,-0.11){8}{\line(1,0){0.36}}
\multiput(65.84,50.88)(0.22,-0.12){12}{\line(1,0){0.22}}
\multiput(68.48,49.45)(0.14,-0.11){17}{\line(1,0){0.14}}
\multiput(70.77,47.52)(0.12,-0.15){16}{\line(0,-1){0.15}}
\multiput(72.64,45.18)(0.11,-0.22){12}{\line(0,-1){0.22}}
\multiput(74.01,42.52)(0.12,-0.41){7}{\line(0,-1){0.41}}
\multiput(74.83,39.63)(0.11,-1.49){2}{\line(0,-1){1.49}}
\multiput(75.06,36.65)(-0.09,-0.74){4}{\line(0,-1){0.74}}
\multiput(74.69,33.67)(-0.12,-0.36){8}{\line(0,-1){0.36}}
\multiput(73.74,30.83)(-0.11,-0.20){13}{\line(0,-1){0.20}}
\multiput(72.25,28.23)(-0.12,-0.13){17}{\line(0,-1){0.13}}
\multiput(70.27,25.98)(-0.15,-0.11){16}{\line(-1,0){0.15}}
\multiput(67.88,24.17)(-0.25,-0.12){11}{\line(-1,0){0.25}}
\multiput(65.19,22.86)(-0.41,-0.11){7}{\line(-1,0){0.41}}
\multiput(62.29,22.12)(-1.50,-0.08){2}{\line(-1,0){1.50}}
\multiput(59.29,21.96)(-0.74,0.11){4}{\line(-1,0){0.74}}
\multiput(56.33,22.39)(-0.31,0.11){9}{\line(-1,0){0.31}}
\multiput(53.51,23.41)(-0.20,0.12){13}{\line(-1,0){0.20}}
\multiput(50.95,24.96)(-0.13,0.12){17}{\line(-1,0){0.13}}
\multiput(48.75,26.99)(-0.12,0.16){15}{\line(0,1){0.16}}
\multiput(46.99,29.42)(-0.11,0.25){11}{\line(0,1){0.25}}
\multiput(45.74,32.15)(-0.11,0.49){6}{\line(0,1){0.49}}
\put(45.07,35.06){\line(0,1){2.99}}
\multiput(44.98,38.06)(0.10,0.59){5}{\line(0,1){0.59}}
\multiput(45.48,41.01)(0.11,0.28){10}{\line(0,1){0.28}}
\multiput(46.57,43.80)(0.12,0.18){14}{\line(0,1){0.18}}
\multiput(48.18,46.33)(0.12,0.12){18}{\line(0,1){0.12}}
\multiput(50.26,48.49)(0.16,0.11){15}{\line(1,0){0.16}}
\multiput(52.73,50.19)(0.28,0.12){10}{\line(1,0){0.28}}
\multiput(55.48,51.36)(0.75,0.12){6}{\line(1,0){0.75}}
\multiput(80.00,49.87)(1.13,-0.11){3}{\line(1,0){1.13}}
\multiput(83.39,49.54)(0.36,-0.11){9}{\line(1,0){0.36}}
\multiput(86.65,48.58)(0.22,-0.11){14}{\line(1,0){0.22}}
\multiput(89.68,47.02)(0.15,-0.12){18}{\line(1,0){0.15}}
\multiput(92.35,44.91)(0.12,-0.14){19}{\line(0,-1){0.14}}
\multiput(94.57,42.34)(0.11,-0.20){15}{\line(0,-1){0.20}}
\multiput(96.27,39.38)(0.11,-0.32){10}{\line(0,-1){0.32}}
\multiput(97.37,36.17)(0.12,-0.84){4}{\line(0,-1){0.84}}
\multiput(97.85,32.80)(-0.09,-1.70){2}{\line(0,-1){1.70}}
\multiput(97.68,29.40)(-0.12,-0.47){7}{\line(0,-1){0.47}}
\multiput(96.86,26.09)(-0.12,-0.26){12}{\line(0,-1){0.26}}
\multiput(95.44,23.00)(-0.12,-0.16){17}{\line(0,-1){0.16}}
\multiput(93.45,20.24)(-0.12,-0.12){20}{\line(-1,0){0.12}}
\multiput(90.97,17.90)(-0.18,-0.11){16}{\line(-1,0){0.18}}
\multiput(88.10,16.08)(-0.29,-0.11){11}{\line(-1,0){0.29}}
\multiput(84.93,14.83)(-0.56,-0.10){6}{\line(-1,0){0.56}}
\put(81.59,14.20){\line(-1,0){3.40}}
\multiput(78.19,14.23)(-0.56,0.11){6}{\line(-1,0){0.56}}
\multiput(74.85,14.89)(-0.29,0.12){11}{\line(-1,0){0.29}}
\multiput(71.70,16.18)(-0.18,0.12){16}{\line(-1,0){0.18}}
\multiput(68.85,18.04)(-0.12,0.12){20}{\line(-1,0){0.12}}
\multiput(66.40,20.41)(-0.11,0.16){17}{\line(0,1){0.16}}
\multiput(64.45,23.20)(-0.12,0.26){12}{\line(0,1){0.26}}
\multiput(63.06,26.31)(-0.11,0.47){7}{\line(0,1){0.47}}
\multiput(62.29,29.62)(-0.06,1.70){2}{\line(0,1){1.70}}
\multiput(62.16,33.02)(0.10,0.67){5}{\line(0,1){0.67}}
\multiput(62.68,36.39)(0.11,0.32){10}{\line(0,1){0.32}}
\multiput(63.83,39.59)(0.12,0.20){15}{\line(0,1){0.20}}
\multiput(65.56,42.52)(0.12,0.13){19}{\line(0,1){0.13}}
\multiput(67.82,45.07)(0.15,0.12){18}{\line(1,0){0.15}}
\multiput(70.51,47.14)(0.23,0.12){13}{\line(1,0){0.23}}
\multiput(73.56,48.67)(0.59,0.11){11}{\line(1,0){0.59}}
\put(40,57){\makebox(0,0)[cc]{$[n]$}}
\put(45.67,23.67){\makebox(0,0)[cc]{$C$}}
\put(50,18.00){\makebox(0,0)[cc]{\small$|C|\!=\!k$}}
\put(101.00,24.00){\makebox(0,0)[cc]{$D$}}
\put(102.00,18.33){\makebox(0,0)[cc]{\small$|D|\!=\!k$}}
\put(67.67,36.00){\makebox(0,0)[cc]{$I$}}
\put(68.33,31.00){\makebox(0,0)[cc]{\small$|I|\!=\!w$}}
\put(71.33,57.3){\makebox(0,0)[cc]{$J$}}
\put(86,57.30){\makebox(0,0)[cc]{\small$|J|\!=\!2k\!-\!2w$}}
\put(70, 2){\makebox(0,0)[cc]{Figure 2.}}
\put(53.67,50.33){\line(-3,-5){8.20}}
\put(59.00,51.67){\line(-3,-5){12.40}}
\put(64.00,51.00){\line(-3,-5){14.80}}
\put(68.67,48.67){\line(-3,-5){15.20}}
\put(62.33,27.67){\line(-3,-5){3.40}}
\put(78.00,49.67){\line(-2,-3){3.78}}
\put(83.67,49.33){\line(-3,-5){8.00}}
\put(89.33,47.00){\line(-3,-5){18.20}}
\put(94.67,42.00){\line(-3,-5){16.60}}
\put(97.67,34.67){\line(-3,-5){11.80}}
\put(69.33,25.00){\line(-2,-3){2.67}}
\put(68.33,57.00){\vector(-1,-1){9.33}}
\put(72.67,55.67){\vector(3,-4){8.00}}
\end{picture}\\

\end{document}